\documentclass[11pt]{amsart}
\usepackage{graphicx,multicol,abraces}
\usepackage{graphics}
\usepackage{amsmath}
\usepackage{amscd}
\usepackage{color}
\usepackage{latexsym}
\usepackage[all]{xy}
\begin{document}
\textwidth 5.5in
\textheight 8.3in
\evensidemargin .75in
\oddsidemargin.75in

\newtheorem{lem}{Lemma}[section]
\newtheorem{claim}[lem]{Claim}
\newtheorem{conj}[lem]{Conjecture}
\newtheorem{defi}[lem]{Definition}
\newtheorem{thm}[lem]{Theorem}
\newtheorem{cor}[lem]{Corollary}
\newtheorem{prob}[lem]{Problem}
\newtheorem{rmk}[lem]{Remark}
\newtheorem{que}[lem]{Question}
\newtheorem{prop}[lem]{Proposition}
\newcommand{\p}[3]{\Phi_{p,#1}^{#2}(#3)}
\def\Z{\mathbb Z}
\def\R{\mathbb R}
\def\g{\overline{g}}
\def\odots{\reflectbox{\text{$\ddots$}}}
\newcommand{\tg}{\overline{g}}
\def\proof{{\bf Proof. }}
\def\ee{\epsilon_1'}
\def\ef{\epsilon_2'}
\title{Homology spheres with $E_8$-fillings and arbitrarily large correction terms}
\author{Motoo Tange}
\thanks{The author was partially supported by JSPS KAKENHI Grant Number 17K14180.}
\subjclass{57R55, 57R65}
\keywords{Definite 4-manifold, $E_8$-bounds, Heegaard Floer correction term, Neumann-Siebenmann invariant}
\address{Institute of Mathematics, University of Tsukuba,
 1-1-1 Tennodai, Tsukuba, Ibaraki 305-8571, Japan}
\email{tange@math.tsukuba.ac.jp}
\date{\today}
\maketitle
\begin{abstract}
In this paper we construct families of homology spheres which bound 4-manifolds with intersection forms isomorphic to $-E_8$.
We show that these families have arbitrary large correction terms.
This result says that among homology spheres, the difference of the maximal rank of minimal sub-lattice of definite filling and the maximal rank of even definite filling is arbitrarily large.
\end{abstract}
\section{Introduction}
\subsection{Definite fillings and homology cobordism invariants.}
If a 3-manifold $Y$ bounds $X$, then we call $X$ {\it filling} of $Y$.
If a filling of $Y$ has a definite (even, or spin) intersection form, then the filling is called {\it definite filling} ({\it even filling} or {\it spin filling} respectively).
Under the assumption that the homology of a filling has no 2-torsions, an even filling is equivalent to spin filling.
If a definite filling has a positive (or negative) definite intersection,
then we call the filling {\it positive-definite filling} (or {\it negative-definite filling} respectively).

Let $Y$ be an integer homology sphere.
Rohlin invariant $\mu(Y)$ is defined to be $\sigma(W)/8\in {\mathbb Z}/2{\mathbb Z}$ for a spin filling $W$ of $Y$.
We can assume that the spin filling $W$ is $H_1(W,{\mathbb Z})=\{0\}$
(we say {\it homologically 1-connected}).
In this article we mainly consider homologically 1-connected definite fillings.

Ozsv\'ath and Szab\'o defined a homology cobordism invariant $d$ in \cite{osabs}.
If a 3-manifold has a negative-definite filling of $Y$, then the $d$-invariant has the following restriction.
\begin{thm}[\cite{osabs}]
Let $Y$ be an integral homology three-sphere, then for each negative-definite four-manifold
$X$ which bounds $Y$, we have the inequality
$$\xi^2+\text{rk}(H_2(X,{\mathbb Z}))\le 4d(Y)$$ 
for each characteristic vector $\xi$.
\end{thm}
Furthermore, if a homology sphere $Y$ has an even negative-definite filling $W$, then $b_2(W)\le 4d(Y)$ is satisfied.
For example $\Sigma(2,3,5)$ is the boundary of the $-E_8$-plumbing.
Here $-E_8$ is the unique unimodular, even, negative-definite, rank 8 quadratic form.
The computation $d(\Sigma(2,3,5))=2$ means that $b_2$ of any even negative-definite filling is at most $8$.
The plumbing realized a filling with $b_2=8$.
If $Y$ has a definite filling with intersection form $nE_8$ for some integer $n$, then 
the filling is called {\it $nE_8$-filling}.

On the other hand, the $d$-invariant of $\Sigma(2,3,7)$ is $0$.
Thus, if there exists an even negative-definite filling, then the $b_2$ has to be $0$.
Hence, the $b_2$ of any positive-definite filling is also zero.
Since $\mu(\Sigma(2,3,7))=1$, it has no homology 1-connected even definite filling.
The plumbing of $\Sigma(2,3,7)$ with all weights $-2$ can give an even filling with even intersection form $(-E_8)\oplus H$, where $H$ is the hyperbolic intersection form.
This filling is a homologically 1-connected even (equivalently spin) indefinite filling.

In \cite{T} the author defined the following invariants.
If $Y$ has an $nE_8$-filling, then we define $g_8$ (or $\underline{g_8}$) to be
$$g_8(Y) = \max\{b_2(W)/8|W:\text{$nE_8$-filling of $Y$},H_1(W)=\{0\}\},$$
$$\underline{g_8}(Y ) = \min\{b_2(W)/8|W:\text{$nE_8$-filling of $Y$},H_1(W)=\{0\}\}.$$
If $Y$ has no $nE_8$-fillings, then $g_8(Y)=-\infty$.

We call the invariant $g_8$ {\it $E_8$-genus}.
If $Y$ has an $nE_8$-filling, then we can immediately see the following bound
\begin{equation}
\label{in}
2g_8(Y)\le |d(Y)|.
\end{equation}
For example, for any integer $n$, $d(\Sigma(2,3,12n+5))=2$ holds.
The author in \cite{T} showed $g_8(\Sigma(2,3,12n+5))=1$ when $0 \le  n \le  13$ or $n = 15$. 
In \cite{T} we gave the examples with $2g_8(Y)=|d(Y)|$.
The simple question is the following:

\begin{que}
Among homology spheres $Y$ with non-negative $E_8$-genus
is $|d(Y)|-2g_8(Y)$ bounded?
\end{que}
We give families of Brieskorn homology spheres to obtain negative answers for this question.
\subsection{Main results}
Here we give the main result:
\begin{thm}
\label{main}
For any integer $n$, Brieskorn homology spheres $\Sigma(|p|,|q|,|r|)$ for a pairwise coprime positive integer triple $(p,q,r)$ below have homologically 1-connected $-E_8$-fillings with $g_8 =-\bar{\mu}= 1$.
{\normalfont \begin{itemize}
\item (i) $(2,8n - 3,14n - 5)$, (ii) $(2,14n + 3,24n + 5) $
\item (iii) $(2,16n + 3,26n + 5)$, (iv) $(2,10n - 3,16n - 5) $
\item (v) $(5,35n - 2,50n - 3)$, (vi) $(5,25n - 2,40n - 3) $
\item (vii) $(3,15n - 2,36n - 5)$, (viii) $(3,9n - 2,24n - 5) $
\item (ix) $(3,21n - 4,36n - 7)$, (x) $(3,27n - 4,48n - 7)$
\item (xi) $(4,28n - 3,64n - 7)$, (xii) $(4,32n - 3,76n - 7)$
\end{itemize}}
\end{thm}
The invariant $\bar{\mu}$ is the Neumann-Siebenmann invariant, which will be defined in Section~\ref{ns}.
These examples can be useful for realizing desired filling restricted by gauge theory.
For example, see recent Scaduto's study \cite{C}.
\begin{thm}
\label{correctiontermbound}
For positive integer $n$ the correction terms of Brieskorn homology spheres (i), (ii) (iii) and (iv) in Theorem~\ref{main} have the following inequalities:
\begin{multicols}{2}
$2\left\lceil\frac{n}{2}\right\rceil \le d(\Sigma(2,8n - 3,14n - 5))$,

$2\left\lceil\frac{n+1}{2}\right\rceil\le d(\Sigma(2,14n+3,24n+5))$,
\end{multicols}
\begin{multicols}{2}
$2\left\lceil \frac{n+1}{2}\right\rceil\le d(\Sigma(2,16n+3,26n+5))$,

$2\left\lceil\frac{n}{2}\right\rceil\le \Sigma(2,10n - 3,16n - 5).$
\end{multicols}
\end{thm}

These theorems say that for any positive integer $n$, the Brieskorn
homology spheres (i), (ii), (iii), and (iv) have $-E_8$-fillings and $d(Y)-2g_8(Y)=d(Y)+2\bar{\mu}(Y)$ are arbitrarily large.
\begin{rmk}
{\normalfont Let $(Y,c)$ be a pair of Seifert rational homology sphere $Y$ and a spin structure $c$.
According to \cite{U2} the $\bar{\mu}(Y,c)$ is the equivalent to the Fukumoto-Furuta invariant $w(Y,c)$.

Manolescu in \cite{Mano} defined homology cobordism invariants $\alpha$, $\beta$, and $\gamma$ in the framework of Pin(2) Seiberg-Witten Floer homology.
A result in \cite{St} says that for any Brieskorn homology sphere $Y$ (with usual orientation)
$\beta(Y)=\gamma(Y)=-\bar{\mu}(Y)$ and $\alpha(Y)=d(Y)/2$ or $d(Y)/2+1$.
Hence, our result means the existence of homology spheres that $\beta(Y)=1$ and $\alpha(Y)$ is arbitrarily large.
}
\end{rmk}
\begin{rmk}
{\normalfont As a conjecture, the inequalities in Theorem~\ref{correctiontermbound} would become the equalities actually for any positive integer $n$.
The evidence is due to Karakurt's program \cite{HFNem}.
Similarly, for positive integer $n$ we predict the following equalities for other Brieskorn homology spheres in Theorem~\ref{main}.
\begin{itemize}
\item For $(p,q,r)=(5,35n-2,50n-3),(5,25n - 2,40n - 3)$,  we have
$$d(\Sigma(p,q,r)) = 6n.$$
\item For $(p,q,r)=(3,15n - 2,36n - 5), (3,9n - 2,24n - 5), (3,21n - 4,36n - 7), (3,27n - 4,48n - 7)$, we have
$$d(\Sigma(p,q,r)) =2n.$$
\item  For $(p,q,r)=(4,28n - 3,64n - 7), (4,32n - 3,76n - 7)$,  we have
$$d(\Sigma(p,q,r))= 4\left(\frac{n}{2}+\left\lceil\frac{n}{2}\right\rceil\right).$$
\end{itemize}
For non-positive integer $n$ we conjecture that for any homology sphere above $d(\Sigma(|p|,|q|,|r|))$  are all $2$. 
}
\end{rmk}
Here we compare the following result by Ue \cite{ue} with the result above.
\begin{thm}[\cite{ue}]
Let $(S,c)$ be a pair of a spherical 3-manifold and a spin structure 
on it.
Then $d(S,c)=-2\bar{\mu}(S,c)$.
\end{thm}
In fact, in \cite{ue} it is shown that the general correction term $d(S,t)$ of spin$^c$ spherical 3-manifold coincides with 
the Fukumoto-Furuta invariant, which is defined by using the index of the Dirac operator of a line bundle over a 4-orbifold.

We remark $\bar{\mu}$ here is defined to be the $\bar{\mu}$ divided by 8 in \cite{ue}.
The examples in Theorem~\ref{main} are homology spheres that $d(\Sigma)+2\bar{\mu}(\Sigma)$ are arbitrarily large.
The relationship between correction term $d$ and Neumann-Siebenmann invariant $\bar{\mu}$ for non-spherical 3-manifolds is known not so many things even in the case of Seifert 3-manifolds. 

\subsection{Other invariants related definite fillings}
\subsubsection{Homologically 1-connected fillings}
In \cite{T} homology cobordism invariant $\frak{ds}(Y)$ is defined to be
\begin{center}
the maximal $b_2(W)/8$ among homologically 1-connected, \\
even definite fillings $W$
\end{center}
of $Y$.
A homologically 1-connected even filling gives a spin filling.
Any invariant related to a kind of filling is defined to be $-\infty$, if there exists no such a kind of fillings.
We define $\frak{o}(Y)$ to be 
\begin{center}
the maximal rank of the minimal definite lattice $L$ that $L\oplus \langle \pm1\rangle^n$ is the intersection form of a homologically 1-connected definite filling $W$
\end{center}
of $Y$.
Here `minimal' means that any square $\pm1$ element is not included in the lattice
and $n$ is some non-negative integer.
By the definition, we have 
\begin{equation}
\label{inf}
8g_8(Y)\le 8\frak{ds}(Y)\le \frak{o}(Y).
\end{equation}
For example, in the case of $Y=\Sigma(2,5,9)$, according to Corollary 1.2 in \cite{C}, 
$g_8(Y)=\frak{ds}(Y)=1$ and $\frak{o}(Y)=12$ holds.
We have the following question.
\begin{que}
Differences of these invariants are bounded or unbounded?
\end{que}
We can show that homology spheres in Theorem~\ref{correctiontermbound} satisfy the following unbounded property.
\begin{cor}
\label{unbounded}
Among homology spheres $Y$, $\frak{o}(Y)-8\frak{ds}(Y)$ is unbounded.
\end{cor}
These integer homology spheres satisfy $\frak{ds}(Y)=g_8(Y)$.
It is not known whether for some integer homology spheres $Y$
the differences $\frak{ds}(Y)-g_8(Y)$ are positive or unbounded.

\subsubsection{General definite fillings}
Define $E(Y)$ to be
\begin{center}
the maximal among $b_2(W)/8$ even definite filling $W$
\end{center}
of $Y$.
Note the filling is possibly non-spin.
In the same way, we define $O(Y)$ to be 
\begin{center}
the maximal rank of minimal sub-lattice of definite fillings $W$
\end{center}
of $Y$ and $G_8(Y)$ to be 
\begin{center}
the maximal $|n|$ among $nE_8$-fillings
\end{center}
of $Y$.
Since fillings for invariants $G_8(Y)$, $E(Y)$, $O(Y)$ do not assume homologically 1-connected, we have $g_8(Y)\le G_8(Y)$, $\frak{ds}(Y)\le E(Y)$, and $\frak{o}(Y)\le O(Y)$.
We have the similar inequality to (\ref{inf}):
\begin{equation}
8G_8(Y)\le 8E(Y)\le O(Y).
\end{equation}

For example, consider the case of $Y=\Sigma(2,3,7)$.
As we mentioned above $Y$ has no even definite fillings with homologically 1-connected,
i.e., $g_8(Y)=\frak{ds}(Y)=-\infty$.
On the other hand, Fintushel and Stern constructed a rational ball that bounds $Y$ in \cite{FS}, 
i.e., $E(Y)=G_8(Y)=0$.
Let $W$ be a negative-definite filling of $Y$
and $\Xi$ the set of characteristic elements in $H_2(W)$.
$d(Y)=0$ means $\underset{c\in \Xi}{\max}(c^2+b_2(W))\le 0$.
The Elkies theorem in \cite{El} concludes the inequality means the equality and the negative-definite lattice must be diagonalized.
For the positive-definite filling of $Y$ one has only to consider negative-definite fillings of $-Y$.
In fact the plumbing lattice of $\Sigma(2,3,7)$ is diagonalized, therefore, $O(Y)=\frak{o}(Y)=0$.


\section*{Acknowledgements}
This study was started by Christopher Scaduto's question in the Gauge Theory  in Fukuoka in 2018 February:
Does $\Sigma(2,5,9)$ bound any 4-manifold with $-E_8$-intersection form?
The author is grateful for motivating to the calculating.
His question is also answered by Scaduto and Golla in \cite{GS}.
The author is grateful to Macro Golla for giving me many useful comments and advice for writing this article.
The last corollary is indicated by  Christopher Scaduto.
\section{Notations and preliminaries}
\subsection{Plumbing diagram}
\label{pdiagram}
For $i=1,2$ let $V_i\to S^2$ be two $D^2$-bundles over $S^2$
or let $V\to S^2$ be a $D^2$-bundle.
For the $D^2$-bundles $V_1$ and $V_2$ we take sub-$D^2$-bundles over each disk in two base spaces $S^2$.
For $V\to S^2$ we take sub-$D^2$-bundles over two disjoint disks in $S^2$.
{\it Plumbing process} is a surgery obtained by identifying two $D^2$-bundles in such a way that one exchanges the roles of their sections and fibers. 
We call the plumbing of $V\to S^2$ {\it self-plumbing}.
Actually, to define the plumbing process we need choose one of the two possibilities of the orientation of the identification as in p.201 in \cite{GS1}.
Since we deal with the tree-type graph only later, then we do not explain the choices.

We define a plumbing diagram (or graph) as explained in \cite{Sav}.
Let $V$ be the set of vertices with a weight function $m:V\to {\mathbb Z}$.
We assign for $v\in V$ the $D^2$-bundle over $S^2$ with the Euler number $m(v)$. 
Let $E$ be the set of edges.
Each edge $\{v,w\}\in E$ of a plumbing diagram means the plumbing process between the $D^2$-bundles over $S^2$.
If $v=w$, then the edge means the self-plumbing.
Hence, if ${\mathbb Z}$-weighted graph $(V,E,m)$ is called {\it plumbing graph} or {\it plumbing diagram}.

Let $(V,E,m)$ be a plumbing diagram.
The plumbing process along a plumbing graph $\Gamma=(V,E,m)$ gives a 4-manifold $P(\Gamma)$
and we call $P(\Gamma)$ a {\it plumbed 4-manifold}. 
The boundary $\partial P(\Gamma)$ of $P(\Gamma)$ is called a {\it plumbed 3-manifold}.
Here $[v]$ is the class represented by the core sphere of the $D^2$-bundle corresponding to the vertex $v$.
The intersection form $(\cdot, \cdot ):H_2(P(\Gamma))\times H_2(P(\Gamma))\to {\mathbb Z}$ on $P(\Gamma)$ is computed from the linear extension of the following definition.
$$([v],[w])=\begin{cases}
m(v)&v=w\\
1&v\neq w\text{ and }\{v,w\}\in E,\\
0&\text{otherwise.}
\end{cases}$$

A tree-type graph with at most one vertex of degree larger than two is called a {\it star-shaped graph}.
A plumbed 3-manifold with a star-shaped graph is called a {\it Seifert
manifold}.
For $i=1,2,\cdots, n$, let $(\alpha^{(i)}_1, \alpha^{(i)}_2,\cdots, a^{(i)}_{r_i})$ be a 
sequence of weights of vertices of the $i$-th branch of a star-shaped graph. We compute the continued fraction for the sequence as follows:
\begin{equation}
\label{conti1}
a_i/b_i = [\alpha_1^{(i)}, \alpha_2^{(i)},\cdots, \alpha_{r_i}^{(i)}],
\end{equation}
where $\alpha^{(i)}_j$ is some integer and $(a_i,b_i)$ are coprime integers.
Here the continued fraction is defined to be
$$[c_1,c_2,\cdots, c_m] =c_1-\frac{1}{c_2-\cdots-\frac{1}{c_m}}.$$
The Seifert manifold with the rational numbers $a_i/b_i (i = 1, 2,\cdots, n)$ and the central weight $e$.
We present such a Seifert manifold as
$$S(e; (a_1, b_1), (a_2, b_2), \cdots, (a_n, b_n)).$$
We call these integers the {\it Seifert invariant}.
Instead of $(a_i,b_i)$, we also present it as $\alpha_1^{(i)} \cdot \alpha_2^{(i)} \cdots \alpha_{r_i}^{(i)}$.
Thus, we denote the plumbing process by dot $\cdot$.
We present several consecutive integers by the power as follows:
$$\cdots  \overbrace{2 \cdot  2 \cdots  2}^{m} \cdots = \cdots  2^{[m]} \cdots.$$
Brieskorn homology sphere $\Sigma(p, q, r)$ is defined to be
$$\{(z_1,z_2,z_3) \in{\mathbb C}^3|z_1^p + z_2^q + z_3^r = 0\} \cap  S^5,$$
where $p, q$, and $r$ are pairwise coprime positive integers.
The manifold is a plumbed 3-manifold with a branch number three star-shaped graph.
The Seifert invariant is
$$S(e; (p, p'), (q, q'), (r, r')),$$
where $e-(p'/p+q'/q+r'/r)=- 1/pqr$.
For example, the plumbing diagram of $\Sigma(2, 3, 5)$ is described as follows:
$$S(-2; (-2)^{[4]}, (-2)^{[2]}, -2).$$

In \cite{T}, the author showed that the Briesrkorn homology sphere $\Sigma(p,q,r)$ whose intersection matrix of minimal plumbed 4-manifold is isomorphic to $-E_8$ is
$\Sigma(2, 3, 5)$ or $\Sigma(3, 4, 7)$.
\subsection{Notations}
\begin{figure}[htbp]
\begin{center}
\includegraphics{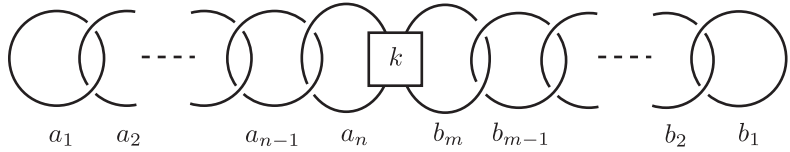}
\caption{A surgery description of (\ref{eq2})}.
\label{sur}
\end{center}
\end{figure}

We explain two new notations as below.
The first notation is the following:
\begin{equation}
\label{eq2}
L(a_1 \cdot a_2 \cdots a_n \cdot^{(k)} b_m \cdot  b_{m-1} \cdots b_1),
\end{equation}
It presents the plumbing as in {\sc Figure}~\ref{sur}.
Here the integers $a_i$,  $b_j$ and $k$ are integers.
A dot with $(k)$ means the two components with framing $a_n$ and $b_m$ 
geometrically link $k$-times. 
Here, the box with the $k$ means the full $k$-twist. 
We call the notation a {\it linear diagram}.

As an example, we consider a linear diagram of $\Sigma(2, 3, 5)$.
Sliding the first branch of
$S(-2; (-2)^{[4]}, (-2)^{[2]}, -2)$ to the third branch, we have the following linear diagram:
\begin{eqnarray*}
S(-2; (-2)^{[4]}, (-2)^{[2]}, -2) &=& L((-2)^{[2]} \cdot  (-2) \cdot  (-2) \cdot^{(-2)} (-4) \cdot  (-2)^{[3]})\\
&= &L((-2)^{[4]} \cdot^{(-2)}(-4) \cdot  (-2)^{[3]}).
\end{eqnarray*}
The 4-manifold having the framed link as in {\sc Figure}~\ref{sur} is called 
{\it 4-manifold having linear diagram (\ref{eq2})}.

\begin{figure}[htbp]
\begin{center}
\includegraphics{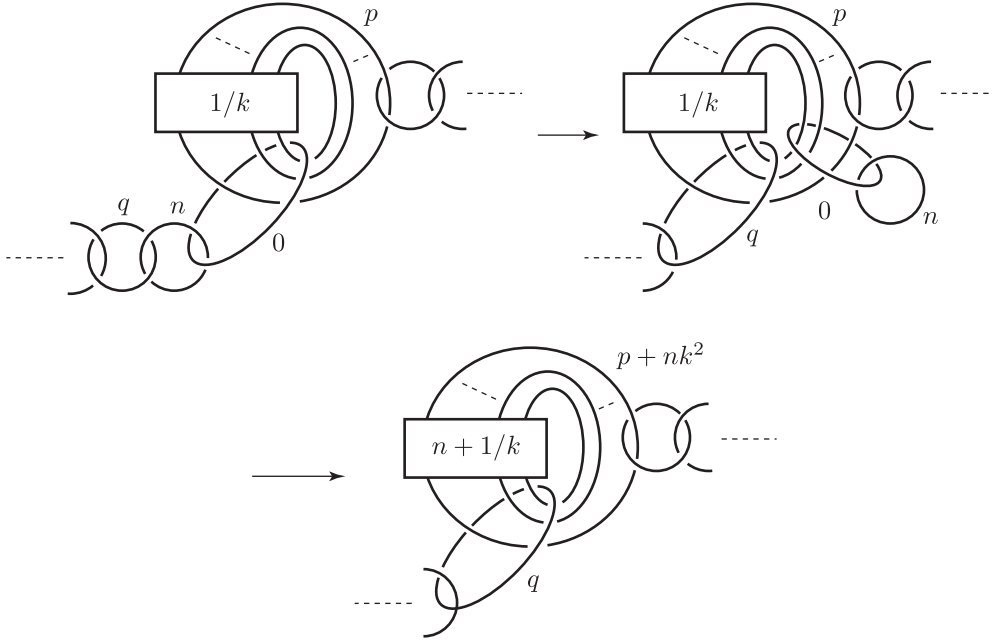}
\caption{A formula of $n$-twisting}
\label{mov}
\end{center}
\end{figure}

The second notation is a linear diagram with torus knot component. 
Consider a linear diagram that the framing of the nearest component to $(k)$ is zero.
Then, the surgery diagram can be deformed as in {\sc Figure}~\ref{mov}.
We present the deformation of the linear diagram as follows:
$$Y=L(\cdots q\cdot n\cdot 0\cdot^{(k)} p\cdots )=L(\cdots q\cdot^{(k)}\underline{p+nk^2}_{(k,nk+1)}\cdots).$$
The component with the underline having index $(k,nk+1)$ stands for the $(k,nk+1)$-torus knot with framing $p+nk^2$.
The last surgery diagram gives a 4-manifold $X$ bounded by $Y$.
Clearly, the intersection form of $X$ is isomorphic to the intersection form of the 4-manifold having the linear diagram as follows:
$$L(\cdots q \cdot^{(k)} (p + nk^2)\cdots).$$
\subsection{An estimate of $\bar{\mu}$-invariant}
\label{ns}
 Neumann-Siebenmann's $\bar{\mu}$-invariant in \cite{NS} is defined for any plumbed 3-manifold $M=\partial P(\Gamma)$.
We set $\Gamma=(V,E,m)$.
We assume that the plumbing graph is tree and $\partial P(\Gamma)$ is a homology sphere.
We define the {\it Wu class} $w(\Gamma) \in  H_2(P (\Gamma), {\mathbb Z})$ as follows:
%
\begin{enumerate}
\item[(1)] The class $w(\Gamma)$ is written by $w(\Gamma)=\sum_{v\in V}\epsilon_v[v]$ for $\epsilon_v=0$ or $1$.
\item[(2)] For any $v\in V$ we have $(w(\Gamma),[v]) \equiv ([v],[v]) \bmod 2$.
\end{enumerate}
Let $\sigma(\Gamma)$ be the signature of the intersection form $(\cdot, \cdot)$ associated with $\Gamma$.
Then we define the $\bar{\mu}$-invariant of $M$ to be
$$\bar{\mu}(M)=\frac{\sigma(\Gamma)-w(\Gamma)^2}{8}.$$
The invariant $\bar{\mu}$ can be extended to any rational plumbed spin 3-manifold $(M,c)$ naturally as in \cite{NS}.
\begin{thm}[\cite{ue}]
\label{uuee}
Suppose that a Seifert rational homology 3-sphere $M$ with spin structure c bounds a negative-definite spin 4-manifold Y with spin structure $c_Y$.
Then
$$b_2(Y )\equiv -8\bar{\mu}(M, c) \bmod 16$$
$$\bar{\mu}(M, c) = \frac{\sigma(\Gamma) - w(\Gamma, c)^2}{8}$$
$$-\frac{8\bar{\mu}(M,c)}{9} \le b_2(Y) \le - 8\bar{\mu}(M,c).$$
\end{thm}
In particular, if a Seifert homology sphere $M$ has a spin negative-definite filling $Y$, then $b_2(Y)\le -8\bar{\mu}(M)$ holds.

\section{The families of Brieskorn homology spheres in Theorem 1.1.}
\subsection{Proof of Theorem \ref{main}}
We prove Theorem\ref{main}.\\
\begin{proof}
The Seifert invariants of
\begin{itemize}
\item $\Sigma(2,14n - 5,8n - 3), \Sigma(2,24n + 5,14n + 3) $
\item $\Sigma(2,26n + 5,16n + 3), \Sigma(2,10n - 3,16n - 5) $
\item $\Sigma(5,35n - 2,50n - 3), \Sigma(5,25n - 2,40n - 3)$ 
\item $\Sigma(3,15n - 2,36n - 5), \Sigma(3,9n - 2,24n - 5) $
\item $\Sigma(2,14n - 5,8n - 3), \Sigma(2,16n - 5,10n - 3) $
\item $\Sigma(3,21n - 4,36n - 7), \Sigma(3,27n - 4,48n - 7)$ 
\item $\Sigma(4,28n - 3,64n - 7), \Sigma(4,32n - 3,76n - 7)$
\end{itemize}
with reverse orientation are as follows:
$S(1;(p_1,1),(p_2,q_2),(p_3,q_3))$.
$$\begin{array}{|c|c|c||c|c|c|}\hline
p_1&(p_2,q_2)&(p_3,q_3)&p_1&(p_2,q_2)&(p_3,q_3)\\\hline
2&(14n-5, 7n-6)&(8n-3, 2)&2& (14n+3, 7n-2)& (24n+5, 6)\\\hline
2&(26n+5, 13n-4)& (16n+3, 4)&2& (10n-3, 5n-4)& (16n-5, 4)\\\hline
5& (35n-2, 28n-3)& (50n-3, 2)&5& (40n-3, 32n-4)& (25n-2, 1)\\\hline
3& (15n-2, 10n-3)& (36n-5, 4)&3& (24n-5, 16n-6)& (9n-2, 1)\\\hline
3& (21n-4, 14n-5)& (36n-7, 4)&3& (48n-7, 32n-10)&(27n-4, 3)\\\hline
4& (28n-3, 21n-4)& (64n-7, 4)&4& (76n-7, 57n-10)& (32n-3, 2)\\\hline
\end{array}$$

We deform the presentation as follows.
The notation $\sim$ stands for a deformation of presentation preserving the intersection form.

\begin{eqnarray*}
&&S(1;2,2\cdot (n+1)\cdot 2^{[6]},(2-4n)\cdot 2)=S(0,2,(-2)\cdot n\cdot 2^{[6]},(2-4n)\cdot 2)\\
&=&L(2^{[6]}\cdot n\cdot 0\cdot^{(2)}(4-4n)\cdot 2)=L(2^{[6]}\cdot^{(2)}\underline{4}_{(2,2n+1)}\cdot 2)\sim L(2^{[6]}\cdot^{(2)}4\cdot 2) 
\end{eqnarray*}
\begin{eqnarray*}
&&S(1;2,2\cdot (n+1)\cdot 4\cdot 2,(-4n)\cdot 2^{[5]})=S(0;2,(-2)\cdot n\cdot 4\cdot 2,(-4n)\cdot 2^{[5]})\\
&=&L(2^{[5]}\cdot (2-4n)\cdot^{(2)}\cdot 0\cdot n\cdot 4\cdot 2)=L(2^{[5]}\cdot \underline{2}_{(2,2n+1)}\cdot^{(2)}4\cdot 2)\\
&\sim& L(2^{[6]}\cdot^{(2)}4\cdot 2)
\end{eqnarray*}
\begin{eqnarray*}
&&S(1;2,2\cdot (n+1)\cdot  4\cdot 2^{[3]},(-4n)\cdot 2^{[3]})=S(0;2,(-2)\cdot n\cdot 4\cdot 2^{[3]},(-4n)\cdot 2^{[3]})\\
&=&L(2^{[3]}\cdot 4\cdot n\cdot 0\cdot^{(2)}(2-4n)\cdot 2^{[3]})=L(2^{[3]}\cdot 4\cdot^{(2)} \underline{2}_{(2,2n+1)}\cdot 2^{[3]})\\
&\sim& L(2^{[3]}\cdot 4\cdot^{(2)}\cdot 2^{[4]})
\end{eqnarray*}
\begin{eqnarray*}
&&S(1;2,2\cdot (n+1)\cdot 2^{[4]},(2-4n)\cdot 2^{[3]})=S(0;2,(-2)\cdot n\cdot 2^{[4]},(2-4n)\cdot 2^{[3]})\\
&=&L(2^{[4]}\cdot  n\cdot 0\cdot^{(2)}(4-4n)\cdot 2^{[3]})=L(2^{[4]}\cdot^{(2)} \underline{4}_{(2,2n+1)}\cdot 2^{[3]})\sim L(2^{[4]}\cdot^{(2)}4\cdot 2^{[3]})
\end{eqnarray*}
\begin{eqnarray*}
&&S(0;-5,5\cdot n\cdot (-7),(2-25n)\cdot 2)=L((-7)\cdot n\cdot 0\cdot^{(-5)}(-3-25n)\cdot 2)\\
&=&L((-7)\cdot^{(-5)}\underline{-3}_{(-5,-5n+1)}\cdot 2)=L(2\cdot 1\cdot (-5)\cdot^{(-5)}\underline{-3}_{(-5,-5n+1)}\cdot 2)
\end{eqnarray*}
Here we slide the $-3$-framed $(-5,-5n+1)$-torus knot component to the component with $(-5)$-framed component.
Then we have a 4-manifold with intersection form of the plumbed 4-manifold for  $S(1;2,2^{[2]},-5)$.
By doing four blow-ups and one blow-down, we have the intersection $E_8$. 

In the same way, the 4-manifolds that the last diagrams in the following equalities
present can be
deformed into 4-manifolds with intersection form $E_8$.
The results of the latter four equalities are a plumbed 4-manifold that reduces to the Seifert invariant $S(2;2^{[2]},2^{[4]},2\cdot 4)=-\Sigma(3,4,7)$.
\begin{eqnarray*}
&&S(0;-5,5\cdot n\cdot (-8),(2-25n))=L((-3-25n)\cdot^{(-5)}0\cdot n\cdot (-5)\cdot 1\cdot 2^{[2]})\\
&=&L(\underline{-3}_{(-5,-5n+1)}\cdot^{(-5)}(-5)\cdot 1\cdot 2^{[2]})
\end{eqnarray*}
\begin{eqnarray*}
&&S(0;-3,3\cdot n\cdot (-5),(2-9n)\cdot 2^{[3]})=L((-5)\cdot n\cdot 0\cdot^{(-3)}(-1-9n)\cdot 2^{[3]})\\
&=&L(2\cdot 1\cdot(-3)\cdot^{(-3)}\underline{(-1)}_{(-3,-3n+1)}\cdot 2^{[3]})
\end{eqnarray*}
\begin{eqnarray*}
&&S(0;-3,3\cdot n\cdot (-8),(2-9n))=L(2^{[4]}\cdot 1\cdot (-3)\cdot n\cdot 0\cdot^{(-3)}(-1-9n))\\
&=&L(2^{[4]}\cdot 1\cdot (-3)\cdot^{(-3)}\underline{(-1)}_{(-3,-3n+1)})
\end{eqnarray*}
\begin{eqnarray*}
&&S(1;2,2\cdot n\cdot (-7),(2-4n)\cdot 2^{[2]})=S(0;-2,2\cdot n\cdot (-7),(2-4n)\cdot 2^{[2]})\\
&=&L(2^{[4]}\cdot 1\cdot (-2)\cdot n\cdot 0\cdot^{(-2)}(2-4n)\cdot 2^{[2]})\\
&=&L(2^{[4]}\cdot 1\cdot (-2)\cdot^{(-2)}\underline{2}_{(-2,-2n+1)}\cdot 2^{[2]})
\end{eqnarray*}
\begin{eqnarray*}
&&S(1;2,2\cdot n\cdot (-6),(2-4n)\cdot 2^{[3]})=S(0;-2,2\cdot n\cdot (-6),(2-4n)\cdot 2^{[3]})\\
&=&L(2^{[3]}\cdot 1\cdot (-2)\cdot n\cdot 0\cdot^{(-2)}(-4n)\cdot 2^{[3]})\\
&=&L(2^{[3]}\cdot 1\cdot (-2)\cdot^{(-2)}\underline{0}_{(-2,-2n+1)}\cdot 2^{[3]})
\end{eqnarray*}
\begin{eqnarray*}
&&S(0;-3,3\cdot n\cdot (-7),(2-9n)\cdot 4)=L(2^{[3]}\cdot 1\cdot(-3)\cdot n\cdot 0\cdot^{(-3)}(-1-9n)\cdot 4)\\
&=&L(2^{[3]}\cdot 1\cdot (-3)\cdot^{(-3)}\underline{(-1)}_{(-3,-3n+1)}\cdot 4)
\end{eqnarray*}
\begin{eqnarray*}
&&S(0;-3,3\cdot n\cdot (-5)\cdot 3,(2-9n)\cdot 2^{[2]})\\
&=&L(2^{[2]}\cdot(-1-9n)\cdot^{(-3)}0\cdot n\cdot (-3)\cdot 1\cdot 2\cdot 4)\\
&=&L(2^{[2]}\cdot\underline{(-1)}_{(-3,-3n+1)}\cdot^{(-3)}(-3)\cdot 1\cdot 2\cdot 4)
\end{eqnarray*}
\begin{eqnarray*}
&&S(0;-4,4\cdot n\cdot (-7),(2-16n)\cdot 4)\\
&=&L(4\cdot (-2-16n)\cdot^{(-4)}0\cdot n\cdot (-4)\cdot 1\cdot 2^{[2]})\\
&=&L(4\cdot \underline{(-2)}_{(-4,-4n+1)}\cdot^{(-4)} (-4)\cdot 1\cdot 2^{[2]})
\end{eqnarray*}
\begin{eqnarray*}
&&S(0;-4,4\cdot n\cdot (-6)\cdot 3,(-2-16n)\cdot 2)\\
&=&L(2\cdot(-2-16n)\cdot^{(-4)}0\cdot n\cdot(-4)\cdot 1\cdot 2\cdot 4)\\
&=&L(2\cdot\underline{(-2)}_{(-4,-4n+1)}\cdot^{(-4)}(-4)\cdot 1\cdot 2\cdot 4)
\end{eqnarray*}
According to the definition of $\bar{\mu}$ as above, computing the $\bar{\mu}$-invariants for these Brieskorn homology spheres, we can see $\bar{\mu}=-1$ easily.
From the description under Theorem~\ref{uuee} we obtain $g_8\le 1$.
Namely, the homology spheres have all $g_8=1$.\\
\hfill$\Box$
\end{proof}

\section{Brieskorn homology spheres with $E_8$-filling and arbitrarily large correction terms.}
\subsection{Heegaard Floer homology and one preparation }
In \cite{osabs} for any spin$^c$ rational homology sphere $(Y,\frak{s})$ the Heegaard Floer homology $HF^+(Y,\frak{s})$ has the following exact sequence:
$$0\to T^+_{d(Y,\frak{s})}\to HF^+(Y,\frak{s})\to HF_{\text{red}}(Y,\frak{s})\to 0.$$
$T^+_s$ is isomorphic to $T^+:={\mathbb F}[U,U^{-1}]/U\cdot {\mathbb F}[U]$ with the minimal degree $s$. 
$HF_{\text{red}}(Y,\frak{s})$ is a finite dimensional torsion ${\mathbb F}[U]$-module.
$d(Y,\frak{s})$ is called the {\it correction term} of $(Y,\frak{s})$.
We call the submodule $T^+_{d(Y,\frak{s})}$ in $HF^+(Y,\frak{s})$ {\it the $T^+$-part} of $HF^+(Y,\frak{s})$.

Here we prepare a lemma to prove Theorem~\ref{correctiontermbound}.
We abbreviate $d(L(p,q),i)$ by $d(p,q,i)$.
Here we define the lens space $L(p,q)$ to be the $p/q$-surgery of the unknot.
The identification of spin$^c$ structures with ${\mathbb Z}/p{\mathbb Z}$ is due to Fig.2 in \cite{osabs}.
Here the $p$-Dehn surgery of a knot $K$ in a homology sphere $Y$ is the surgery 
$(Y\setminus S^1\times D^2)\cup V$, where $V\cong S^1\times D^2$.
Here the attaching meridian of the new solid torus $V$ is mapped to $p\cdot [m]+[l]\in H_1(\partial(S^1\times D^2))$ 
where $m$ is the meridian of $K$ and $l$ is the homologically trivial longitude of $K$.
We denote the $p$-Dehn surgery of a knot in $\Sigma$ by $\Sigma(p)$.
\begin{lem}
\label{onepre}
Let $\Sigma$ be a homology sphere and $K\subset \Sigma$ a knot.
For some positive integer $p$, if $\Sigma(p-1)$ is an L-space and $\Sigma(p)$ is a lens space $L(p,q)$, then the correction term $d(\Sigma)$ is computed as follows.
\begin{equation}
\label{eq1}
d(\Sigma)=\max\left\{d\left(p,q,ki+c\right)-d(p,1,i)|0\le i<p\right\},
\end{equation}
where $k$ is the dual class of $[\tilde{K}]\in H_1(L(p,q),{\mathbb Z})$
and $c=(k+1+p)(k-1)/2$,
where $\tilde{K}$ is the surgery dual of the lens space surgery.
\end{lem}
Let $C$ be a core circle of the genus one Heegaard decomposition.
Then the {\it dual class} $k$ is defined $k[C]=[\tilde{K}]\in H_1(L(p,q),{\mathbb Z})$.
The dual class is used in the situation of the integral lens space surgery of a homology sphere in \cite{T3}.

\begin{proof}
We use the following surgery exact sequence (Corollary 9.13 in \cite{osabs}):
$$\cdots \to HF^+(\Sigma)\to  HF^+(\Sigma(p-1)) \overset{G^+}{\to} HF^+(\Sigma(p))\overset{F^+}{\to} HF^+(\Sigma)\to  \cdots.$$
Since $\Sigma(p-1)$ and $\Sigma(p)$ are L-spaces and
the corresponding map $F^\infty$ on $HF^\infty$ is surjective, 
$F^+$ is also surjective onto the  $T^+$-part in $HF^+(\Sigma)$.
The map $F^+$ is induced from the cobordism $\Sigma(p)$ to 
$\Sigma$
obtained by attaching a 0-framed 2-handle along the meridian of $\tilde{K}$.
The spin$^c$ structures on $\Sigma(p)$ are identified with ${\mathbb Z}/p{\mathbb Z}$ due to the description in p.213 in \cite{osabs}.
For any integer $j$ with $0\le j<p$ consider the surgery exact sequence in Theorem 9.19 in \cite{OS2}:
$$\cdots \to HF^+(\Sigma)\to  HF^+(\Sigma(0),[j]) \to HF^+(\Sigma(p),j)\overset{F_{j}^+}{\to} HF^+(\Sigma)\to  \cdots.$$
$F_{j}^+$ is a component of $F^+$ restricted to the spin$^c$ structure $j$.
It is also a sum of homogeneous maps $f_{i}^+$ with respect to the spin$^c$  cobordism from $(\Sigma(p),j)$ to the unique spin$^c$ manifold on $\Sigma$.
Namely, $F_j^+$ is described by the sum $F_{j}^+=\sum_{j\equiv i\bmod p}f^+_i$.
The degree shift of $f_i^+$ is $(4p-(2i-p)^2)/(4p)$ due to \cite{osabs}.
The maximal degree shift among $\{f_i^+|j\equiv i\bmod p\}$ is $(4p-(2j-p)^2)/(4p)=-d(p,1,j)$.
Since $F_{j}^+$ is a surjective $U$-equivariant map, for $0\le j<p$ we have
$$d(\Sigma)\ge d(p,q,kj+c)-d(p,1,j).$$
The 1-1 correspondence ${\mathbb Z}/p{\mathbb Z}\to \text{Spin}^c(L(p,q))$ in
Corollary 7.5 in \cite{osabs} is described by $ki+c$.
See \cite{T2}.

Suppose that $d(\Sigma)>d(p,q,kj+c)-d(p,1,j)$ for any integer $j$ with $0\le j<p$.
Then any element with the minimal degree in $HF^+(\Sigma(p))$ is included in the kernel of $F^+=\sum_{0\le j<p}F_j^+$.
Thus the kernel of $F^+$ includes at least $p$ components.
On the other hand, for a sufficient large number $N$,
$\text{ker}(F^+)/(U^N=0)$ is ($p-1$)-fold direct sum of $T^+$ from the exact sequence of the version of $HF^\infty$.
Hence, this implies that in the image of $G^+$ there is a torsion ${\mathbb F}[U]$-module by at least one component.
However, since $\Sigma(p-1)$ is an L-space, the image of $G^+$ has no torsion ${\mathbb F}[U]$-module.
This is a contradiction.
Therefore for some $j$, $d(\Sigma)=d(p,q,kj+c)-d(p,1,j)$ holds.
\hfill$\Box$
\end{proof}

\subsection{The $d$-invariants for the four families of Brieskorn homology spheres.}
We prove Theorem~\ref{correctiontermbound}.

\begin{proof}
The Seifert presentations of Brieskorn homology spheres from (i) to (iv) in Theorem~\ref{main} are the below:
$$
\begin{array}{|l|l|}\hline
\text{(i)}&S(1;2,2\cdot (-n+1)\cdot 7,(4n-1)\cdot 2)\\\hline
\text{(ii)}&S(1;2,2\cdot (-n)\cdot (-3)\cdot 2,(4n+1)\cdot 6)\\\hline
\text{(iii)}&S(1;2,2\cdot (-n)\cdot (-3)\cdot 4,(4n+1)\cdot 4)\\\hline
\text{(iv)}&S(1;2,2\cdot (-n+1)\cdot 5,(4n-1)\cdot 4)\\\hline
\end{array}
$$
Let $\Sigma_n$ be one of Brieskorn homology spheres parametrized by $n$ in the list above.
We do 0-surgery and +1-surgery of the homology sphere along the meridian of the singular fiber of multiplicity 2.
We call the meridian $K_n$.
Note the coefficients $0$ and $1$ are the framing of the unknot $K_n$ in the diagram.
The $0$-surgery and $1$-surgery give lens spaces $L(r_n,s_n)$ and $L(p_n,q_n)$.
The results are the lens spaces in the list below.
$$
\begin{array}{|l|l|l|}\hline
&\text{0-surgery}\ (r_n,s_n)&\text{1-surgery}\ (p_n,q_n)\\\hline
\text{(i)}&(56n^2-41n+7,8n^2-7n+2)&(56n^2-41n+8,8n^2-7n+1)\\\hline
\text{(ii)}&(168n^2+71n+7,72n^2+27n+4)&(168n^2+71n+8,72n^2+27n+1)\\\hline
\text{(iii)}&(208n^2+79n+7,48n^2+17n+2)&(208n^2+79n+8,48n^2+17n+1)\\\hline
\text{(iv)}&(80n^2-49n+7,16n^2-13n+4)&(80n^2-49n+8,16n^2-13n+1)\\\hline
\end{array}
$$
These examples satisfy $p_n=r_n+1$.
As a result, the $0$-surgery means a positive $r_n$-Dehn surgery along $K_n$.

%
%
%
%

We set $d_n:=d(\Sigma_n)$.
Here using Lemma~\ref{onepre}, we compute the lower bound of $d_n$.
We argue the case of (i) only.
Other cases are able to be proven by similar arguments.
Let $\Sigma_n$ be the Brieskorn homology sphere in the type (i).
We set $p_n=56n^2-41n+8$, $q_n=8n^2-7n+1$, $k_n=14n-5$ and $c_n\equiv 42n^2-29n+4\bmod p_n$.
The $k_n$ is the dual class in the lens space $L(p_n,q_n)$ which presents $K_n$.

Here we set $i=\lfloor\frac{q_n+1}{2}\rfloor-n$.
Then by modulo $p_n$ we have
$$k_ni+c_n\equiv\begin{cases}\frac{7n-5}{2}&n\text{: odd}\\-\frac{7n}{2}&n\text{: even.}\end{cases}$$

If $n$ is an odd number, then by using the reciprocity formula in \cite{osabs}, 
we have 
$$d\left(L(p_n,q_n),\frac{7n-5}{2}\right)=\frac{224n^3+8n^2-95n+25}{4p_n}$$
and
$$d(L(p_n,1),i)=-\frac{52n^2-37n+7}{4p_n}$$
Thus we have
$$d(L(p_n,q_n),k_ni+c_n)-d(L(p_n,1),i)=n+1.$$

If $n$ is an even number, then we have 
$$d\left(L(p_n,q_n),-\frac{7n}{2}\right)=\frac{224n^3-216n^2+73n-8}{4p_n}$$
$$d(L(p_n,1),i)=-\frac{52n^2-41n+8}{4p_n}.$$
Thus we have
$$d(L(p_n,q_n),k_ni+c_n)-d(L(p_n,1),i)=n.$$

Therefore we have $d_n\ge 2\lceil \frac{n}{2}\rceil$.
\hfill$\Box$
\end{proof}
\subsection{Proof of Corollary~\ref{unbounded}.}
Here we prove Corollary~\ref{unbounded}.

\begin{proof}
Let $\{\Sigma_n\}$ be one family of homology spheres in Theorem~\ref{correctiontermbound}.
Since the correction term is positive, $\Sigma_n$ is no even positive-definite filling of $\Sigma_n$.
A Seifert homology sphere $\Sigma_n$ 
has a negative-definite plumbing $\Sigma_n=\partial P(\Gamma_n)$.
Due to Corollary 1.5 in \cite{OS3}, $4d(\Sigma_n)=\underset{c\in \Xi}{\max}(c^2+\text{rank}(\Gamma_n))$ holds, where $\Xi_n$ is the set of characteristic classes in $\Gamma_n$.
If for a non-negative integer $N$, $\Gamma_n=\Gamma_n'\oplus\langle -1\rangle^N$ and $\Gamma_n'$
is a minimal sub-lattice, then $\underset{c\in \Xi}{\max}(c^2+\text{rank}(\Gamma_n))$ is decomposed into the sum of the two maximal values
according to the direct sum.
Thus, we have 
$$\underset{c\in \Xi}{\max}(c^2+\text{rank}(\Gamma_n))=\underset{c\in \Xi'}{\max}(c^2+\text{rank}(\Gamma_n'))\le \text{rank}(\Gamma_n'),$$
where $\Xi'$ is the set of characteristic elements of $L'$.
Since $d(\Sigma_n)$ has no upper bound,
the maximal of the rank of minimal lattice $\Gamma_n'$ is also unbounded.
On the other hands, $\bar{\mu}(Y_n)=1$ and our construction of filling of $\Sigma_n$ implies that $\frak{ds}(Y_n)=1$.
Therefore $\frak{o}(\Sigma_n)-\frak{ds}(Y_n)$ is unbounded.
\hfill$\Box$
\end{proof}

%
%
%

\noindent
Motoo Tange\\
University of Tsukuba, \\
Ibaraki 305-8502, Japan. \\
tange@math.tsukuba.ac.jp

\end{document}